\documentclass[12pt]{article}

\usepackage[letterpaper, left=2.5cm, top=2cm, right=2.5cm, bottom=2.5cm ]{geometry}
\usepackage{authblk}
\usepackage[T1]{fontenc}

\usepackage[bookmarks=closed]{hyperref}
\usepackage{breakurl}

\usepackage{amsmath,amsfonts,amssymb}
\usepackage{bm}
\usepackage[caption=false]{subfig}
\usepackage{textcomp}
\usepackage{stfloats}
\usepackage{hyperref}
\usepackage{url}
\usepackage{verbatim}
\usepackage{siunitx}
\usepackage{graphicx}
\usepackage{amssymb, amsthm}
\usepackage{cite}

\usepackage{enumitem}
\usepackage{multirow}

\usepackage{todonotes}

\usepackage{siunitx}

\usepackage{mathtools}
\usepackage{booktabs}

\allowdisplaybreaks

\title{Heat diffusion blurs photothermal images with increasing depth}

\date{}

\author{Peter Burgholzer}

\affil{Research Center for Non Destructive Testing (RECENDT)\authorcr
4040 Linz, Austria\authorcr
E-mail:  \texttt{peter.burgholzer@recendt.at}
 }

\author{G\"unther Mayr}
\author{Gregor  Thummerer}

\affil{Josef Ressel Center for Thermal NDE of Composites\authorcr
University of Applied Sciences Upper Austria\authorcr
4600 Wels, Austria}

\author{Markus Haltmeier}

\affil{Department of Mathematics, University of Innsbruck\authorcr
6020 Innsbruck, Austria}

\begin{document}

\maketitle

\begin{abstract}

In this tutorial, we aim to directly recreate some of our ``aha'' moments when exploring the impact of heat diffusion on the spatial resolution limit of photothermal imaging. Our objective is also to communicate how this physical limit can nevertheless be overcome and include some concrete technological applications. Describing diffusion as a random walk, one insight is that such a stochastic process involves not only a Gaussian spread of the mean values in space, with the variance proportional to the diffusion time, but also temporal and spatial fluctuations around these mean values. All these fluctuations strongly influence the image reconstruction immediately after the short heating pulse. The Gaussian spread of the mean values in space increases the entropy, while the fluctuations lead to a loss of information that blurs the reconstruction of the initial temperature distribution and can be described mathematically by a spatial convolution with a Gaussian thermal point-spread-function (PSF). The information loss turns out to be equal to the mean entropy increase and limits the spatial resolution proportional to the depth of the imaged subsurface structures. This principal resolution limit can only be overcome by including additional information such as sparsity or positivity. Prior information can be also included by using a deep neural network with a finite degrees of freedom and trained on a specific class of image examples for image reconstruction.

\end{abstract}

\section{\label{sec:Introduction}Introduction}
Photothermal imaging, and also photoacoustic imaging, involves (scattered)  light causing a rise in the temperature of sub-surface light-absorbing structures, such as blood vessels in human tissue or carbon fibers in an epoxy material. In photothermal imaging, the internally absorbed heat diffuses from the absorbing structures to the sample surface, where it can be measured as a change in surface temperature. In photoacoustic imaging, the temperature rise causes a sudden increase in pressure, which produces an acoustic wave that can be detected at the surface with an ultrasound transducer. Both surface signals, the detected temperature increase and the pressure signal, can be used to reconstruct an image of the internal light absorbing structures in a non-destructive way. The achievable spatial resolution in imaging is directly related to the information content of the thermal diffusion wave or the acoustic wave\cite{JAPtutorial}. For linear image reconstruction methods, the spatial resolution is the width of the point-spread-function (PSF), which is the reconstructed image of a point-like heat or pressure source.\par

Unfortunately, heat diffusion from the internal structure to the sample surface ``annihilates'' a lot of information, which, as we have shown earlier, corresponds to entropy production during heat diffusion \cite{Burgholzer.2013,Burgholzer.2015}. In thermal imaging, this leads to a spatial resolution proportional to the depth of a point heat source below the surface. In photoacoustic imaging, entropy production due to acoustic attenuation is lower and leads for liquids, where the acoustic attenuation is proportional to the square of the acoustic frequency, to a degradation of spatial resolution proportional to the square root of depth \cite{Patch.2007}. Therefore, for a higher imaging depth, acoustic imaging shows a significantly better spatial resolution than thermographic imaging.\par

\begin{figure}
    \centering 
    \includegraphics[width=0.5\textwidth]{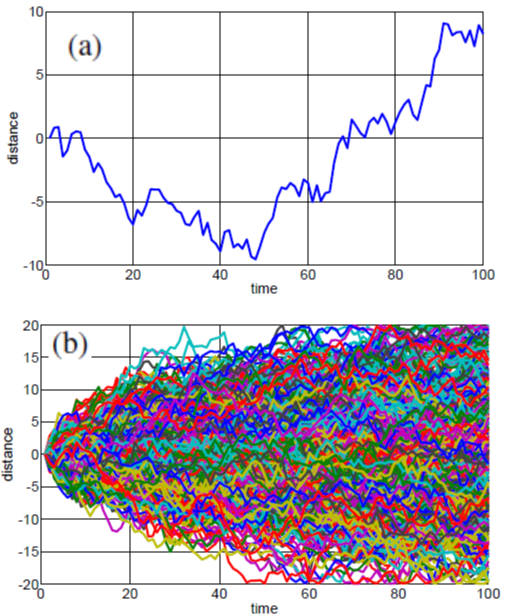}
    \caption{One realization of the Wiener process with the particle starting at $x=0$ (a) and 500 realizations using reflecting boundaries at a distance of $\pm$20 (b).}
    \label{fig:Fig_random_walk}
\end{figure}

In this tutorial, we first describe the one-dimensional diffusion process as a random walk, which can be described mathematically as simple stochastic process. Of course, this is a rough approximation to the real physical processes of heat diffusion, where phonons can carry a whole distribution of energy, but two essential properties of thermal diffusion waves can already be recognized in this random walk model: The mean values of the occupation numbers spread in space, which can be described by the diffusion equation, and the actual occupation numbers fluctuate around these mean values. This fluctuation is described as the variance of the stochastic process. The diffusion process spreads the particles or the heat in space and at the same time induces temporal and spatial fluctuations around the mean values of the Gaussian spreading, which are the reason for an information loss being equal to an increase of entropy. For a macroscopic object, the temperature fluctuations are very small, but for thermographic reconstruction, the time reversal of heat diffusion turns out to be an highly ill-posed inverse problem that amplifies these fluctuations exponentially with increasing time. This limits the spatial resolution proportional to the depth of the subsurface structures, and this principal resolution limit can only be overcome by including additional information such as sparsity or positivity \cite{Thummerer.2020b}. Prior information can be also included by using a trained deep neural network for image reconstruction, either directly from the measured data (end-to-end reconstruction by deep learning) or by calculating first the virtual wave, which corresponds to  the same initial heat or pressure distribution as the measured surface temperature, but is a solution of the wave equation (hybrid reconstruction by the virtual wave concept and deep learning)\cite{SPM}.\par

\begin{figure}[htb!]
    \centering
    \includegraphics[height=0.62\columnwidth]{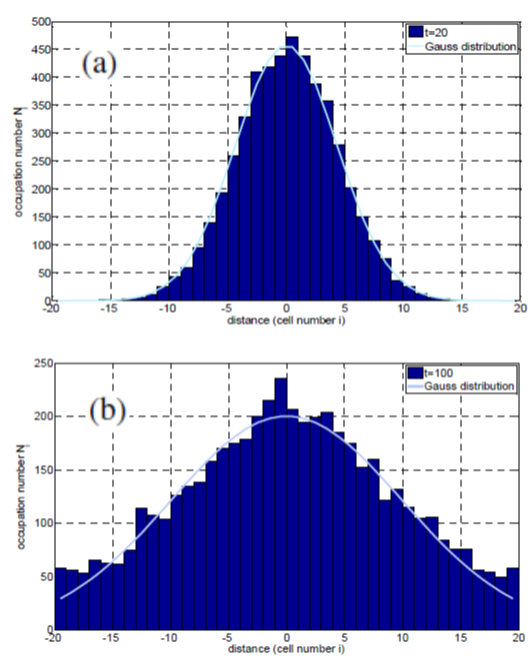}
    \includegraphics[height=0.62\columnwidth]{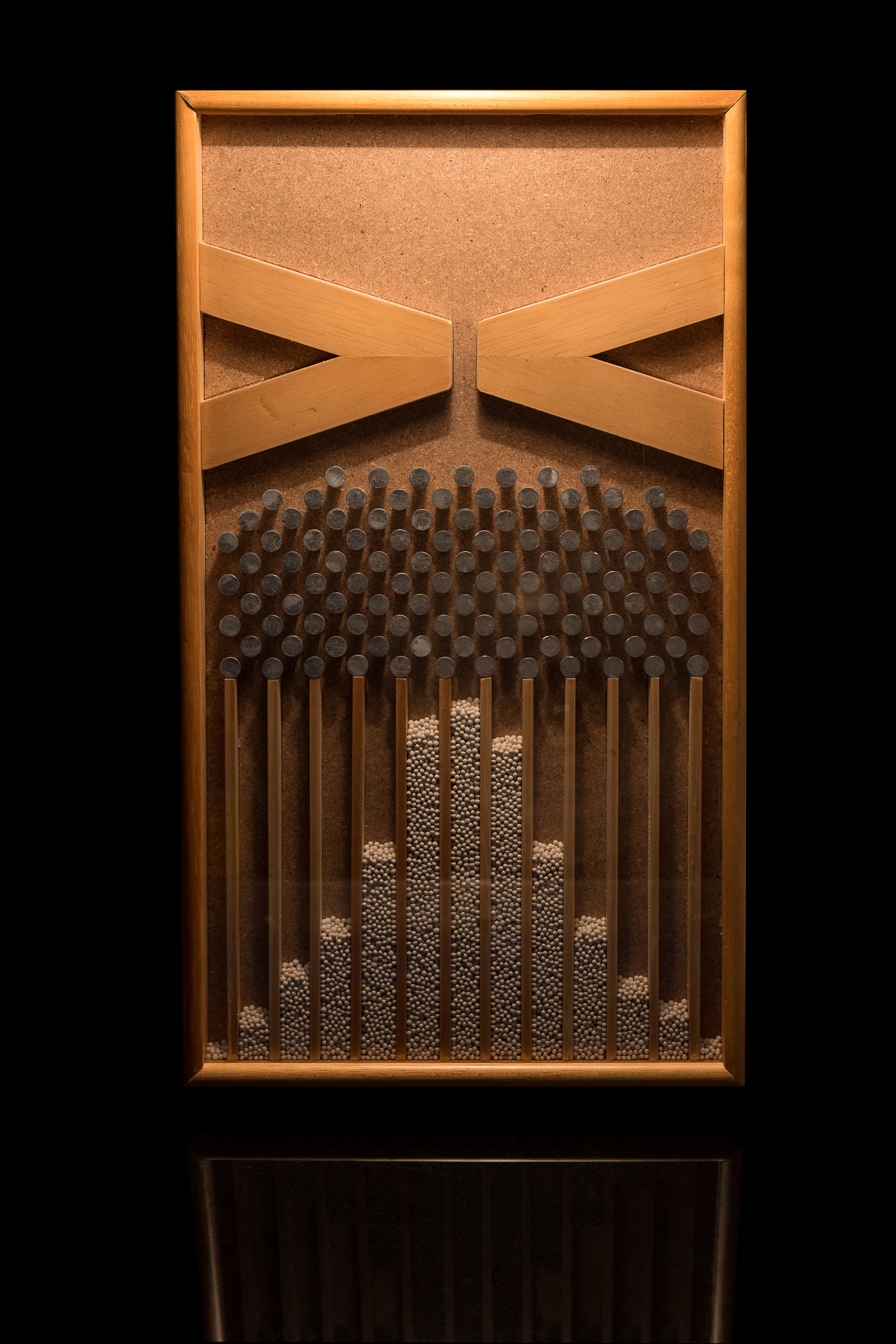}
    \caption{Left: Occupation number representation of the Wiener process for 5000 particles starting at $x=0$ at $t=20$ (a) and $t=100$ (b). For comparison the Gaussian distribution with a variance proportional to time is shown. For $t=100$ the reflecting boundaries at a distance of $\pm$20 can be already recognized by occupation numbers significantly increasing compared to the Gaussian distribution. Right: Galton board to demonstrate the occupation number representation. Licensed as Creative Commons BY-SA 4.0: Matemateca (IME/USP)/Rodrigo Tetsuo Argenton.}
    \label{fig:Fig_GaussGalton}
\end{figure}

\section{\label{sec:random walk}Heat diffusion described as a random walk}
In one dimension, the traditional way to define a random walk is to allow the walker to take steps at each time interval at which he must step either backward or forward, with equal probability (Fig. \ref{fig:Fig_random_walk}).  Described mathematically as a stochastic process, for the so-called Wiener process, the distribution density of the walking distance is a Gaussian distribution with the mean value constant in time and equal to the starting position, and a variance proportional to time $t$ (e.g. Gardiner \cite{Gardiner.1985}). Therefore, an initially sharp distribution spreads out with time (Fig. \ref{fig:Fig_random_walk}(b) and Fig. \ref{fig:Fig_GaussGalton} left).  In addition to these simulated results, such a stochastic process can be physically realized using a Galton board (Fig. \ref{fig:Fig_GaussGalton}). The Galton board consists of a vertical board with nested rows of pins. Beads are dropped from the top and bounce either to the left or to the right when they hit the pins. They are collected into bins at the bottom, where the height of bead columns accumulated in the bins approximate a Gaussian bell curve. For the following mathematical description of the one-dimensional diffusion process, the individual bins are the cells numbered such that the central cell has $i=0$, and the occupation number $N_i(t)$ is the number of beads collected in a given pin. The time $t$ can be identified with the number of rows of pins of the Galton board.

We assume that the diffusion of the particles is independent from each other, with $p_i(t)$ giving the probability that a certain particle is in cell $i$ at time $t$ when starting from cell $i=0$ at $t=0$. $X_{ij}(t)$ is a new random variable, defined to be one if particle $j$ is in cell $i$ and zero in all other cases. For every particle, $X_{ij}(t)$ is equal to one with probability $p_i(t)$, independent from particle number $j$. Therefore, the mean value and variance of this random variable are the same for all the particles:
\begin{align*}
       \overline{X_{ij}}(t) 
       &= \sum_{i} X_{ij}(t)p_i(t)=p_i(t) \\
       \operatorname{Var}(X_{ij}(t))&=\overline{X_{ij}^2}-(\overline{X_{ij}})^2 \\  &  = \sum_{i} X_{ij}^2 p_i - p_i^2  =p_i(t)-p_i(t)^2 \,.
\end{align*}

To derive these equations, it was used that $X_{ij}^2 = X_{ij}$ equals one if the particle $j$ is in the cell $i$, and equals zero otherwise.
With the occupation numbers $N_i=\sum_{j=1}^N X_{ij}$ the mean value and variance follows from independence of the particle number $j$:
\begin{align}
       \overline{N_i}(t) &= N p_i(t);
       \label{Eq:MeanOccupation} \\[0.2em]
       \operatorname{Var}(N_i(t))&=N (p_i(t)-p_i(t)^2)\nonumber \\ &=\overline{N_i}(t) (1-p_i(t)) \nonumber \\  &\approx \overline{N_i}(t) \; \text{for} \; p_i(t) \ll 1 \,.
        \label{Eq:Occupation}
\end{align}

\begin{figure}
    \centering
    \includegraphics[width=0.7\columnwidth]{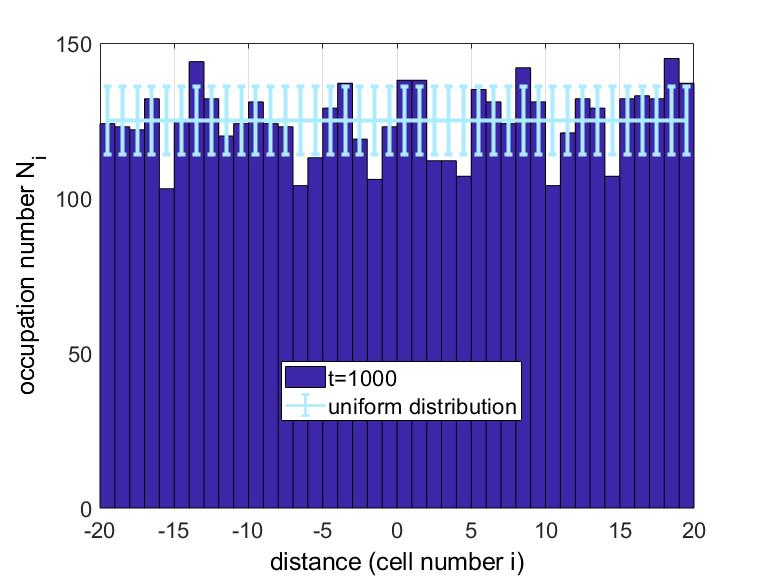}
    \caption{Occupation number representation of the Wiener process for $N=5000$ particles starting at $x=0$ at $t=1000$ using 40 cells. At that time the influence of the reflecting boundaries at a distance of $\pm$20 dominates. For comparison the uniform distribution with the mean value 125 (5000 particles divided by 40 cells) and the erroroverline at $\pm$ the standard deviation is shown. The standard deviation is the square root of the variance given in Eq. (\ref{Eq:Occupation}) using $p_i=1/40$ for the uniform distribution.}
    \label{fig:Occupation}
\end{figure}

A simple model for one dimensional heat diffusion using adiabatic boundary conditions is a one dimensional random walk with reflecting boundaries. Such adiabatic boundary conditions can be realized for a thermally isolated sample, for which the normal derivative of the temperature vanishes at all times at the sample boundaries. The walking particles could be thought as the phonons in the sample volume. Giving them all the same energy is a rough approximation, but the essential details of the influence of fluctuations and the connection of entropy production and information loss can be seen already for such a model system.

\begin{figure}
    \centering
    \includegraphics[width=0.7\columnwidth]{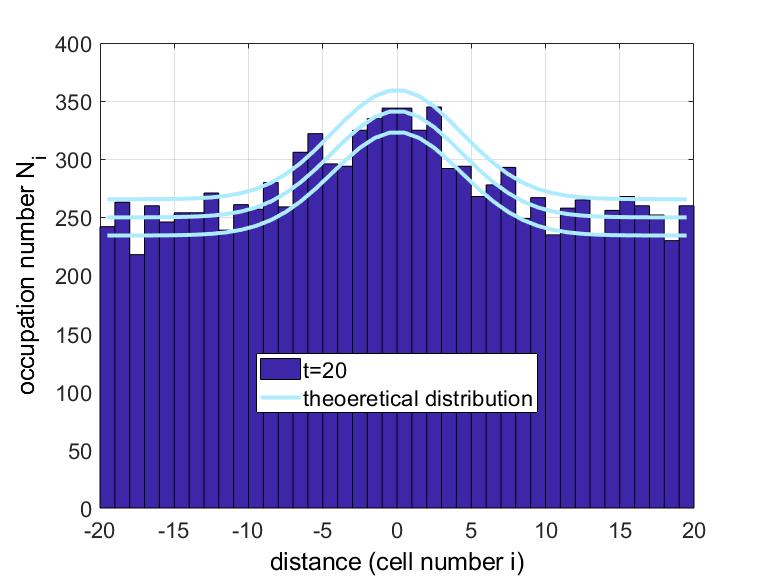}
    \caption{Occupation number representation of the Wiener process for for $N_\text{equi} =10 000, N_\text{cell} =40$ and $N_0=1000$ at a time $t=20$. For comparison the mean value ± the standard deviation (square root of variance) from Eq. (\ref{Eq:Wiener}) is shown, where the variance in a good approximation is constant at $N_\text{equi}/N_\text{cell}$.}
    \label{fig:Wiener}
\end{figure}

Starting from an equilibrium modeled by a uniform distribution (Fig. \ref{fig:Occupation}) with $N_\text{equi}$ particles in $N_\text{cell}$ cells, at the time zero, $N_0$ particles are added at cell $i=0$. Compared to Eq.~ (\ref{Eq:Occupation}) we have now two different distributions: the equilibrium distribution, which is constant in time and space (for every cell $i$), and the Gaussian distribution for the added particles (at least up to a time when the reflecting boundaries can be still neglected). As the movement of the individual particles is assumed to be independent one gets for the mean value and variance:
\begin{align}
   \overline{N_i}(t) &= \frac{N_\text{equi}}{N_\text{cell}} + N_0 p^\text{Gauss}_i(t) \\[0.2em] \nonumber
       \operatorname{Var}(N_i(t))&=\overline{N_i}(t)\bigg(1-\frac{1}{N_\text{equi}+N_0}\overline{N_i}(t)\bigg) 
       \\  \qquad &\approx \overline{N_i}(t) \approx \frac{N_\text{equi}}{N_\text{cell}} \,.
        \label{Eq:Wiener}       
\end{align}
\noindent
The last approximation is valid if the equilibrium occupation number $N_\text{equi}/N_\text{cell}$ is much higher than the mean diffusion term $N_0 p^\text{Gauss}_i(t)$ and the variance $\operatorname{Var}$ gets constant in space and time. An example for $N_\text{equi} =10 000, N_\text{cell} =40$ and $N_0=1000$ is shown in Fig. \ref{fig:Wiener}.

Despite the independence of the movement of the individual particles, the occupation numbers are not uncorrelated. One example is that if a cell has a high occupation number, then the occupation number of the neighboring cells increases shortly afterwards. The covariance of the occupation number in cell $i$ at time $t$ and cell $j$ at a later time $t+\tau$ is:
\begin{align}
       \operatorname{Cov}(N_i(t),&N_j(t+\tau)) \nonumber \\[0.2em]
       &=\sum_{k=1}^N \operatorname{Cov}(X_{ik}(t),X_{jk}(t+\tau)) \nonumber  \\[0.2em] \nonumber
       &= \sum_{k=1}^N \overline{X_{ik}(t) X_{jk}(t+\tau)}-N p_i(t)p_j(t+\tau)\\[0.2em] 
       &= N p_i(t) (p_{i-j}(\tau)-p_j(t+\tau)) \,.
       \label{Eq:covariance}
\end{align}
For Eq. (\ref{Eq:covariance}) it was used that the probability for a particle to be in cell $j$ at time $t+\tau$ when it was in cell $i$ at time $t$ is $p_{i-j}(\tau)$.

The time development of the occupation numbers $N_i$ is a Gauss-Markov process and can be described by a Langevin equation (see e.g. Groot and Mazur\cite{Groot.1984}), which we have used in previous publications to describe the information loss for a kicked Brownian particle\cite{JAPtutorial} or the diffusion of heat\cite{Burgholzer.2013}. During the infinitesimal time $\text{d}t$ the infinitesimal change of the occupation number in the cell $i$ is $\text{d}N_i$:
\begin{equation}
       \text{d}N_i=\frac{1}{2}(N_{i-1}+N_{i+1}-2N_i)\text{d}t+\sqrt{\operatorname{Var} \: \text{d}t}\epsilon_t(0,1) \,,
\end{equation}
where $\epsilon_t(0,1)$ is a random variable with standard distribution (mean value zero and variance one), which describes uncorrelated white noise. The reflecting boundaries are represented by:
\begin{align}
       \text{d}N_1&=\frac{1}{2}(N_2-N_1)\text{d}t+\sqrt{\operatorname{Var} \: \text{d}t}\epsilon_t(0,1) \\
       \text{d}N_{N_\text{cell}}&=\frac{1}{2}(N_{N_\text{cell}-1}-N_{N_\text{cell}})\text{d}t+\sqrt{\operatorname{Var} \: \text{d}t}\epsilon_t(0,1) \,.
\end{align}
All the occupation numbers can be combined in an $N_\text{cell}$ dimensional vector \textbf{N}:
\begin{equation}
       \text{d}\textbf{N}=-\textbf{M} \cdot \textbf{N} \text{d}t+\sqrt{\operatorname{Var} \: \text{d}t} \bm{\epsilon}_t(0,1) \,.
    \label{Eq:Langevin}  
\end{equation}
\noindent
The matrix \textbf{M} in Eq. (\ref{Eq:Langevin}) is a square matrix with the dimension equal to the number of cells:
\begin{equation}
\textbf{M}=-\frac{1}{2}\left( \begin{array}{rrrrrr}
-1 & 1 & 0 &\cdots & 0 & 0\\
1 & -2 & 1 &\ddots & 0 & 0\\
0 & 1 & -2 & 1 &\ddots & \vdots\\
\vdots & \ddots &\ddots &\ddots &\ddots &0\\
0 & \ddots & \ddots & 1 & -2 & 1 \\
0 & 0 & \cdots & 0 & 1 & -1\\
\end{array}\right) 
\end{equation}

\noindent
Eq. (\ref{Eq:Langevin}) can be seen also as a discretized version (finite differences) of the one-dimensional diffusion equation
\begin{equation}
       \frac{\partial N(x,t)}{\partial t} = \alpha \frac{\partial^2 N(x,t)}{\partial x^2}
        \label{Eq:Diffusion}
\end{equation}
with a unit spacing between the discrete cells, a diffusion coefficient of $\alpha=\frac{1}{2}$, adiabatic boundary conditions, and a white noise generating term. The Fourier transform in space ``k-spac'') of the diffusion equation yields
\begin{multline}
       \frac{\partial \hat N (k,t)}{\partial t} = -\alpha k^2 \hat N (k,t) \; \text{with} \; \hat N (k,t)\\=\int N(x,t) e^{ikx}\,\text{d}x \,,
\end{multline}
which solves the diffusion equation (\ref{Eq:Diffusion}) for the  initial condition $N(x,t=~0)=N_0(x)$
\begin{equation}
       N(x,t) = \frac{1}{2 \pi}\int \hat N_0(k) e^{-ikx} e^{-k^2 \alpha t}\,\text{d}k.
       \label{Eq:Fourier}
\end{equation}

To solve the Langevin equation (\ref{Eq:Langevin}), the matrix \textbf{M} can be diagonalized by an adequate coordinate transform, which turns out to be the Fourier transform in space and which is the same as for the solution of the diffusion equation (\ref{Eq:Diffusion}). For the discretized random walk with reflecting boundaries this transformation is in fact the discrete cosine transform $F=\text{dct}(\text{diag}(1,1,...1))$ because of the adiabatic boundary conditions. $F \cdot M \cdot F^t$ with $F^t$ being the transposed matrix of $F$ results in a diagonal matrix with the singular values $\gamma_k$ in the matrix diagonal:
\begin{multline}
       \gamma_k=2 \text{sin}^2(k/2) \approx k^2/2 \; \\\text{with} \; k=\frac{\pi}{N_\text{cell}} (0,1,...,N_\text{cell}-1) \,.
       \label{Eq:SingularValues}
\end{multline}
The last approximation is for a continuous k-space when the cell volume gets zero and $N_\text{cell}$ gets to infinity and is identical to $\alpha k^2$ with $\alpha=1/2$ for the one-dimensional diffusion equation in k-space. The singular values are proportional to the square of the wavenumber $k$ and are the time constants for the exponential decay of $\hat N(k,t) = \hat N_0(k) e^{-k^2 \alpha t}$ in time (compare to Eq. (\ref{Eq:Fourier})). Therefore variations with a higher wavenumber (lower wavelength) decay faster. Wavenumber $k=0$ shows no decay and represents the equilibrium temperature, which is reached after a long time. The fluctuations and their Variance $\operatorname{Var}(\cdot)$ is the same in k-space and in real space, because the noise energy must be the same in both domains (Parseval's theorem). In section \ref{sec:ResolutionLimit}, we show that the spatial resolution limit at a certain time $t$ is given by the wavenumber $k_\text{cut}$ due to information loss by fluctuations, where $\hat N (k_\text{cut},t)$ is equal to the noise level given by the square root of the variance $\operatorname{Var}$. For wavenumbers with $k>k_\text{cut}$ the $\hat N(k,t)$ is below the noise level and $\hat N_0(k) = \hat N(k,t) e^{+k^2 \alpha t}$ cannot be reconstructed any more.\par

\section{\label{sec:Laplace}How probability comes into play}

In the previous section, diffusion was introduced by a stochastic process such as random walk, where in each time interval the next step is chosen with a certain probability. But how does probability come into play for diffusion? In an earlier publication, we have used statistical physics and non-equilibrium thermodynamics to describe the heat diffusion, entropy, and information during photothermal imaging \cite{JAPtutorial}. One main result was that the information content about the initial state just after the excitation pulse in the diffusive system decreases due to fluctuations. The loss of information in the diffusive system is equal to the increase of entropy, called entropy production, with the entropy expressed by the Shannon or Gibbs formula
\begin{equation}
       S(t) = -k_B \sum_{i} p_i(t) \ln(p_i(t)),
       \label{Eq:entropy}
\end{equation}
where the mean occupation number is $\overline{N_i}(t) = N p_i(t)$ in cell $i$ at time $t$ (Eq. (\ref{Eq:MeanOccupation})) and $k_B$ is the Boltzmann constant. For real heat diffusion, the actual distribution with all its fluctuations can hardly be determined, but for obtaining the information loss, the average entropy production can be used, which for heat diffusion in macroscopic samples is, to a very good approximation, the heat diffused to the surroundings divided by temperature\cite{JAPtutorial}.\par

The temporal evolution of the particles is described microscopically by mechanical equations, which are deterministic. But how does probability come into play in particle diffusion, for example? Here we are confronted by the "demon" first articulated by Pierre Simon Laplace in 1814, which is a hypothetical observer that knows the position and momentum of every molecule, and together with the deterministic time evolution this would be also true for the past and the future. Fortunately, this "demon can be exorcised" by moving to a quantum perspective on statistical mechanics, as described in an excellent overview in Physics Today by Katie Robertson on "demons in thermodynamics" \cite{PhysicsToday}. For Laplace's demon the entropy $S$ in Eq. (\ref{Eq:entropy})) would be zero for all times, because the state of the system is known with certainty. No probability distributions as described in thermodynamics would be necessary. But in the quantum case, probabilities are already an inherent part of the theory, so there is no need to add ignorance to the picture. In other words, the probabilities from statistical mechanics and quantum mechanics turn out to be one and the same\cite{PhysicsToday}.\par

\begin{figure}
	\centering
    \includegraphics[width=0.6\columnwidth]{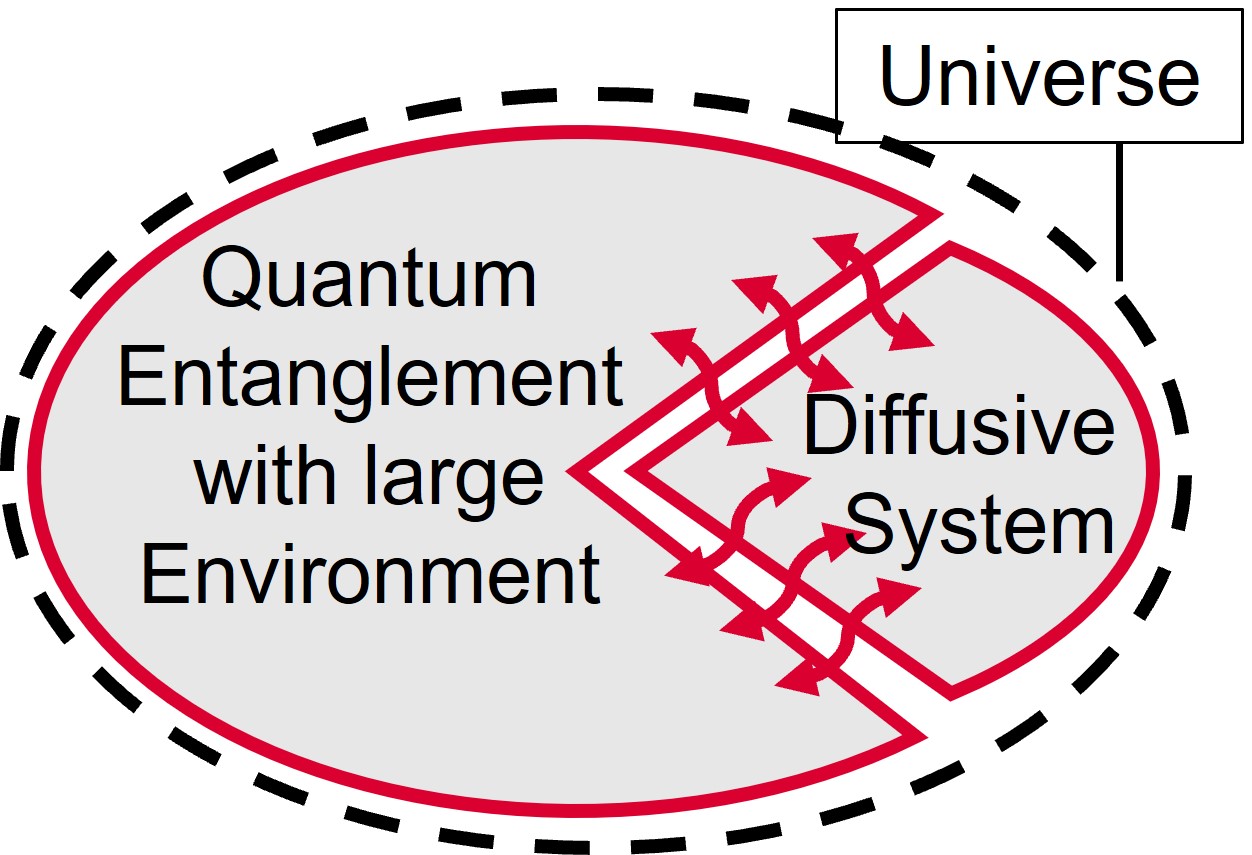}
    \caption{Quantum entanglement explains how probability comes into play. A global state of a large isolated system, called universe, which is in a quantum pure state, can be separated into a small diffusive system entangled with a larger environment. The reduced state of the diffusive system is given by tracing out the environment, and is almost indistinguishable from the thermodynamic state given by statistical mechanics.}
    \label{fig:entanglement}
\end{figure}

Three publications by different authors in 2006 \cite{PhysicsToday1, PhysicsToday2, PhysicsToday3} describe that quantum entanglement of the diffusive system with the environment is the key, and no "ignorance" probabilities are necessary in the description of the diffusion process. They considered the global state of a large isolated system, called the "universe", to be a quantum pure state. Therefore, there is no lack of knowledge about this state, and the entropy is zero. However, if only a small part of the universe is considered, which we call the " diffusive system ", its state will not remain pure due to quantum entanglement with the rest of the universe, which is called the " environment "(Fig. \ref{fig:entanglement}) \cite{PhysicsToday2}. If one of the two systems is taken on its own, it will be in an intrinsically uncertain state known as mixed state. Assuming that the surrounding environment, e.g. the laboratory at room temperature, is sufficiently large, then for almost any pure state that the composite system is in, the diffusive system will be in a state very close to the state it would be assigned by traditional statistical mechanics, and which was described by the stochastic random-walk process in the previous section. The quantum description thus leads to a probability distribution indistinguishable from that of statistical mechanics \cite{PhysicsToday}. Probabilities are an intrinsic and inescapable part of quantum mechanics. When it describes the diffusive system taken on its own, Laplace's demon cannot know more than anyone else. According to Fig. \ref{fig:entanglement}, the probabilities describing the uncertainties (or ignorance) in the diffusion process result from the fact that the information "transported" to the environment by entanglement is lost.\par

\section{\label{sec:ResolutionLimit}Information loss limits the spatial resolution of the reconstruction}

\begin{figure}
\centering
    \includegraphics[width=0.6\columnwidth]{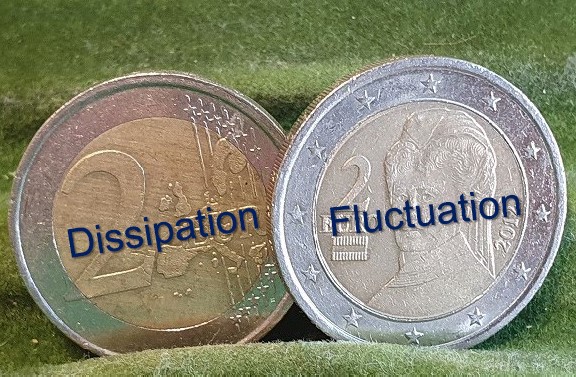}
    \caption{The diffusion process generates spreading in space (dissipation) and fluctuations simultaneously, which can  be  thought  as  two  sides  of  the  same  coin  and  are connected by the fluctuation – dissipation theorem.}
    \label{fig:coins}
\end{figure}

In section \ref{sec:random walk} it was shown that a stochastic process such as random walk leads to spreading in space (dissipation), which is accompanied by an increase in entropy, and at the same time to fluctuations. Fluctuations and mean entropy production can be thought as two sides of the same coin and are connected by the fluctuation – dissipation theorem  (Fig. \ref{fig:coins}). For systems near thermal equilibrium in the linear regime such relations between entropy production and fluctuation properties have been found by Callen \cite{Callen.1952}, Welton \cite{Callen.1951}, and Greene \cite{Greene.1952}. This fluctuation-dissipation theorem is a generalization of the famous Johnson \cite{Johnson.1928} - Nyquist \cite{Nyquist.1928} formula in the theory of electric noise. It is based on the fact that in the linear regime the fluctuations decay according to the same law as a deviation from equilibrium following an external perturbation. We could show from non-equilibrium thermodynamics that the mean entropy production due to the spreading of heat or particles (dissipation) is equal to the information loss due to fluctuations \cite{JAPtutorial}. This connects dissipation, e.g. described by the singular values $\gamma_k$ in Eq. (\ref{Eq:SingularValues}), with the fluctuations given by the variance $\operatorname{Var}$. Eq. (\ref{Eq:Langevin}) reads in k-space for all wavenumbers $k$
\begin{equation}
       \text{d}\hat N_k=- \gamma_k \hat N_k \text{d}t+\sqrt{\sigma_k^2 \: \text{d}t}\epsilon_t(0,1),
       \label{Eq:Langevin_Fourier}
\end{equation}
where the amplitude of the noise for each wavenumber $k$ is given by $\sigma_k^2=2 \gamma_k Var$. This fluctuation-dissipation relation connects directly the exponential decay of the $\hat N(k,t) = \hat N_0(k) e^{-\gamma_k t}$, which shows the dissipative behaviour, with the fluctuations represented by the noise term $\sigma_k^2$ in the Langevin equation (\ref{Eq:Langevin_Fourier}).\par

The diffusion process increases the entropy, but also leads to fluctuations that reduce the information content and thus the spatial resolution for a subsequent reconstruction process. For thermal diffusion imaging, we assume that immediately after the excitation pulse, the absorbed heat is present only at a specific location, e.g., in cell zero. It diffuses with time ($t=20$ (a) and $t=100$ (b) according to Fig. \ref{fig:Fig_GaussGalton}), and because of the fluctuations, reconstruction of the initial distribution from the distribution at a later time is possible only with increasing uncertainty. At time $t=1000$, equilibrium is reached (Fig. \ref{fig:Occupation}) and all information about the initial location is lost.\par

\begin{figure}
    \centering
    \includegraphics[width=0.6\columnwidth]{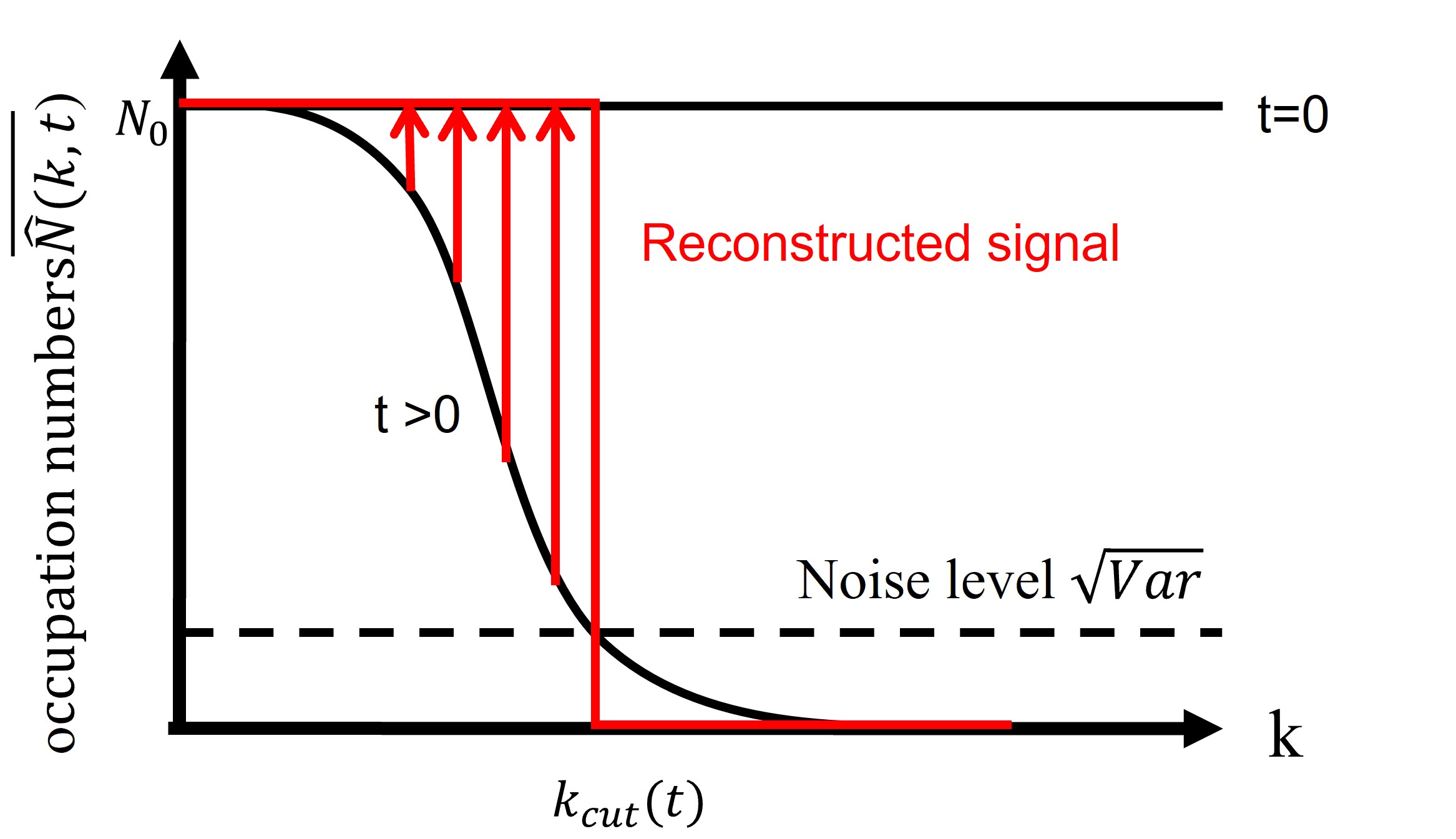}
    \caption{With an initially fully localized distribution $N_0(x) = N_0 \delta (x)$ at $t=0$, which is a horizontal line at $N_0$ in k-space, Eq. (\ref{Eq:Langevin_Fourier}) gives for each wavenumber $k$ an exponential decay of the mean values. After a certain time $t$ only the wavenumbers $k<k_\text{cut}$ where $\overline {\hat N(k,t)} = N_0 e^{-k^2 \alpha t} > \sqrt{\operatorname{Var}}$ can be reconstructed, which gives a rectangular function with $N_0$ from $k=0$ to $k=k_\text{cut}$ and zero above, as shown in red .(Adapted from P.Burgholzer, Int. J. Thermophys 36, 2328–2341, 2015; licensed under Creative Commons Attribution (CC BY) license.)}
    \label{fig:k-space}
\end{figure}
The information loss and the degradation of the spatial resolution for the initial distribution just after the excitation pulse can be quantified in k-space. With an initially fully localized distribution $N_0(x) = N_0 \delta (x)$, Eq. (\ref{Eq:Langevin_Fourier}) gives for each wavenumber $k$ an exponential decay of the mean values. After a certain time $t$ only the wavenumbers $k<k_{cut}$ where $\overline {\hat N(k,t)} = N_0 e^{-k^2 \alpha t} > \sqrt{\operatorname{Var}}$ can be distinguished from noise (Fig. \ref{fig:k-space}) and can be used for the reconstruction \cite{Burgholzer.2015}. The reconstructed signal $N_r(x,t)$ is the inverse Fourier transform of this rectangular function which gives a sinc-function,
\begin{equation}
       N_r(x,t) = \frac{1}{2 \pi}\int_{-k_\text{cut}}^{+k_\text{cut}} N_0 e^{-ikx} \text{d}k = \frac{N_0}{ \pi} \frac{\text{sin} (k_\text{cut}(t) x)}{ x} \,.
       \label{Eq:Fourier_initial}
\end{equation}
With the signal-to-noise ratio $SNR_k=N_0/ \sqrt{\operatorname{Var}}$ in k-space we get
\begin{equation}
      k_\text{cut}= \sqrt{ \frac{\operatorname{ln}(SNR_k)}{\alpha t}} \,,
\end{equation}
with the natural logarithm $\operatorname{ln}(\cdot)$. The spatial resolution for the reconstruction of the initial distribution from a distribution at the time $t$ is the "width" of the reconstructed signal and is taken as half of the "wavelength" at the wavenumber $k_\text{cut}$
\begin{equation}
      \delta_r(t) = \frac{\pi}{k_\text{cut}(t)} = \pi \sqrt{ \frac{\alpha t}{\ln(SNR_k)}}.
      \label{Eq:delta_resolution}
\end{equation}

\section{\label{sec:ThermographicImaging}Thermal point-spread-function}

In general, photothermal imaging does not measure the temperature at all locations $T(x)$ of the sample at a given time~$t$, as assumed in the previous section, but at a specific location $x=0$, which is usually the sample surface, the temperature is measured for the entire time $T(x=0,t)$ \cite{JAPtutorial}. Instead of a Fourier transform to k-space the one-dimensional heat diffusion equation
\begin{equation}
       \frac{\partial T(x,t)}{\partial t} = \alpha \frac{\partial^2 T(x,t)}{\partial x^2}
        \label{Eq:HeatDiffusion}
\end{equation}
with the initial temperature distribution $T_0(x) \equiv T(x,t=0)$ is solved in the $\omega$ - space not by assuming a real wavenumber and a complex frequency as for solving Eq. (\ref{Eq:Diffusion}), but for real frequency $\omega$ and a complex wavenumber $\sigma(\omega)$, and by inserting the "thermal wave" as an Ansatz into the heat diffusion Eq.~(\ref{Eq:HeatDiffusion})
\begin{equation}
    T_{\omega}(x,t) = \Re \left(T_0 \exp \left(i \left(\sigma(\omega)x-\omega t\right)\right)\right) ,
\end{equation}
\noindent
where $\Re$ is the real part, $\sigma(\omega) = \sqrt{i \omega / \alpha } \equiv (1+i)/\mu $ is the complex wave number, and $T_0$ is the amplitude, we get
\begin{equation}
    T_{\omega}(x,t) = T_0 \exp \left(-\frac{x}{\mu}\right)\cos \left( \frac{x}{\mu} - \omega t \right), 
\end{equation}	
\noindent	
with $\mu(\omega) \equiv \sqrt{2 \alpha / \omega} $ is defined as the thermal diffusion length \cite{Salazar.2006}. The thermal wave is strongly damped because the real and imaginary parts of the complex wavenumber $\sigma(\omega)$ are equal and the amplitude is reduced by a factor of $1/e$ after propagation of the thermal diffusion length. The wavenumber or spatial frequency is $k(\omega) \equiv 1 / \mu(\omega) = \sqrt{\omega / 2 \alpha}$. Similar as in k-space, for  frequencies larger than the truncation frequency $\omega_\mathrm{cut}$ the amplitude of this wave components are damped below the noise level, and with the truncation frequency $\omega_\mathrm{cut}$, the truncation diffusion length $\mu_\mathrm{cut}$, and the signal-to-noise ratio $SNR$ in $\omega$ - space we get
\begin{eqnarray}
    {SNR} \; && \exp \left(-\frac{x}{\mu_\mathrm{cut}}\right) = 1, \, \text{or} \nonumber \\ \omega_\mathrm{cut} &&= 2 \alpha \left(\frac{\ln({SNR})}{x}\right)^2.
    \label{Eq:kcut}
\end{eqnarray}

Similar to Eq. (\ref{Eq:delta_resolution}) the spatial resolution 
\begin{equation}
\delta_r(x) = \pi \sqrt{ \frac{2 \alpha} {\omega_\mathrm{cut}}} = \frac{\pi x}{\ln({SNR})}
    \label{Eq:delta_resolution_thermal}    
\end{equation}
\noindent
is half of the wavelength of the truncation frequency. For thermographic reconstruction, the resolution limit increases linear with depth $x$. In comparison, $\delta_r(t)$ increases with the square root of $t$ (Eq. (\ref{Eq:delta_resolution})).\par

In the frequency domain, a short-time excitation is the Fourier transform of a delta pulse where all frequencies have the same amplitude. The thermal point-spread function (PSF) is the inverse Fourier transform, where the frequency integral is taken up to $\omega_\text{cut}$, which assumes that below $\omega_\text{cut}$ the full signal can be reconstructed and for frequency components above the cutoff frequency no signal information is available \cite{Burgholzer.2017b}. Therefore, the signal-to-noise ratio $SNR$, or more precisely, due to the exponential decrease of the amplitude, the natural logarithm $\ln$ of the $SNR$, plays an important role for the spatial resolution. The width of the one-dimensional thermal PSF is proportional to $x/\ln(SNR)$, as derived in Eq. (\ref{Eq:delta_resolution_thermal}). 
For frequency components lower than $\omega_\text{cut}$, the signal amplitude at the sample surface is less than the noise. In deriving spatial resolution, it is assumed that such frequency components cannot be detected. Frequency components above the noise level are used for reconstruction. From information theory and non-equilibrium thermodynamics we could show that the information content for frequencies larger than $\omega_\text{cut}$ is so small that their distribution cannot be distinguished from the equilibrium distribution \cite{JAPtutorial}. If the mean energy and thus the entropy is proportional to the square of the signal amplitude, the cutoff frequency $\omega_\text{cut}$ from this information criterion is always the same frequency at which the signal amplitude becomes smaller than the noise amplitude.\par

One-dimensional signals can be realized by layered structures, and there the noise can easily be reduced by increasing the area covered by the detector, which corresponds to an enhanced averaging. For higher dimensions, such as image reconstruction in three dimensions, this is no longer possible because the signal is different for different detection points on the surface. There is a trade-off between small measurement pixels for better spatial resolution and larger measurement pixels with less noise and better $SNR$, which also affects the spatial resolution of the reconstruction.\par

\begin{figure}[htb]
	\centering
    \includegraphics[width=0.6\columnwidth]{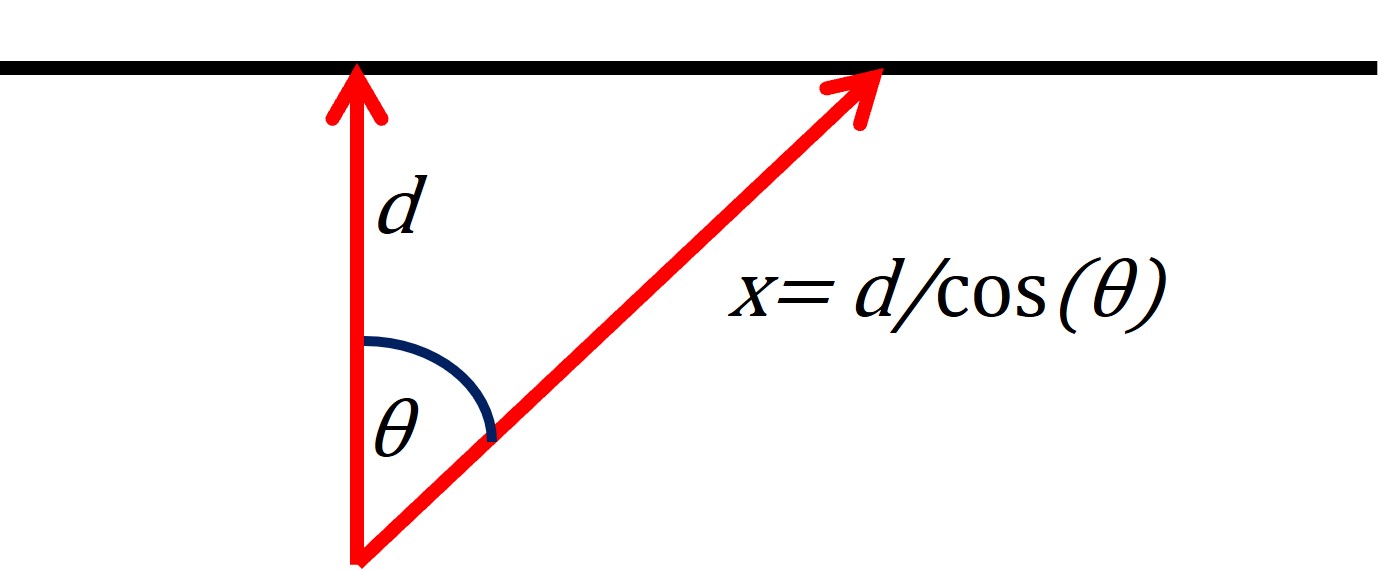}
    \caption{Point source at a depth $d$ below the surface plane, The length $x$ for the thermal signal to reach the surface plane depends on the angle $\theta$. (Adapted from P. Burgholzer et al., Appl. Phys. Lett. 111, 031908, 2017; licensed under Creative Commons Attribution (CC BY) license.)}
    \label{fig:PSF1}
\end{figure}

\begin{figure}[htb]
\centering
    \includegraphics[width=0.5\columnwidth]{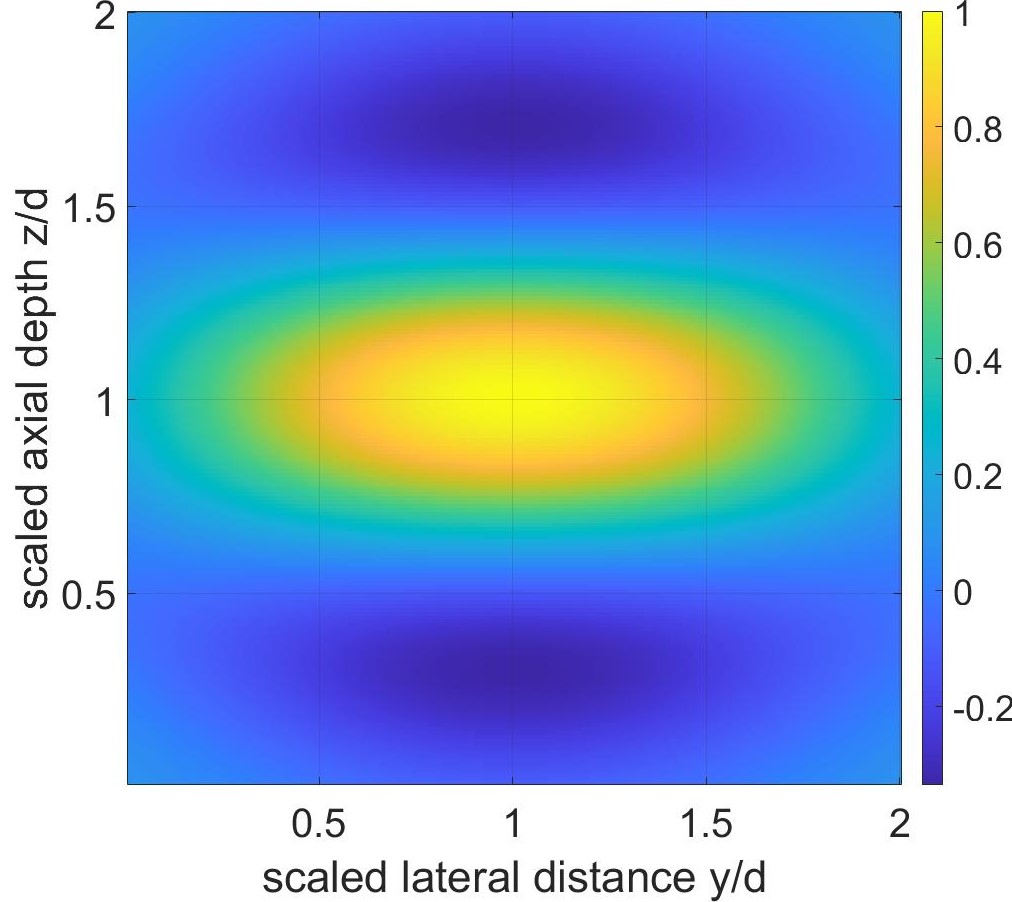}
    \caption{Two-dimensional (or cross section of the three-dimensional) PSF in real space: The lateral resolution (horizontal direction) is 2.4 times the axial resolution (vertical direction). The axial resolution is the same as in the one-dimensional case given in Eq. (\ref{Eq:delta_resolution_thermal}). (Adapted from P. Burgholzer et al., Appl. Phys. Lett. 111, 031908, 2017; licensed under Creative Commons Attribution (CC BY) license.)}
    \label{fig:PSF2}
\end{figure}

To obtain the two- or three-dimensional thermal PSF, we assume a point source to be embedded in a homogeneous sample at a depth $d$ and a planar detection surface. The distance to the surface $x=d/\cos(\theta)$ depends on the angle $\theta$ (Fig. \ref{fig:PSF1}). Using Eq. (\ref{Eq:kcut}) the wavenumber $k_\text{cut}= \sqrt{\omega_\text{cut} / 2 \alpha}$ depends on the direction $\theta$:
\begin{equation}
k_\text{cut}= \sqrt{\frac{\omega_\text{cut}} {2 \alpha}}= \frac{\ln({SNR}) \cos(\theta)}{d}.
    \label{Eq:2d}    
\end{equation}
Inverse Fourier transformation back to real space gives the thermal PSF, shown for $SNR=1000$ in Fig. \ref{fig:PSF2}. It shows the principle resolution limit for a certain $SNR$, like the Abbe limit in optics for a certain wavelength. The axial depth resolution is limited by $k_\text{cut}$ (Eq. (\ref{Eq:2d})) at $\theta =0$, which is the same as in the one-dimensional case given in Eq. (\ref{Eq:delta_resolution_thermal}), and is $\pi d / \ln({SNR})$, and ranges in the scaled depth units in Fig. \ref{fig:PSF2} from $z/d=0.7726$ till $z/d=1.2274$ for a $SNR=1000$. For a $SNR=100$ this would range from 0.66 till 1.34.\par

\section{\label{sec:Overcome}Overcoming the resolution limit}

The resolution limit $\delta_r$ given in Eq. (\ref{Eq:delta_resolution_thermal}) degrades proportional to the imaging depth $x$. One way to improve the resolution when imaging a structure at a certain depth is to increase the $SNR$, which could be achieved by increasing the power of the excitation laser. This is only possible to a certain extent to prevent the sample from being destroyed by too much heating. As mentioned above, for layered structures for which the one-dimensional thermal diffusion equation applies, the $SNR$ can be improved by measuring and averaging over a larger area. For two-dimensional structures, such as parallel cylindrical objects, it is still possible to perform averaging in the direction of the parallel axis and thereby increase the $SNR$. \par
\begin{figure}[h!]
\centering
\includegraphics[trim={8.0cm 6.5cm 8cm 2.3cm},clip,width=0.9\textwidth]{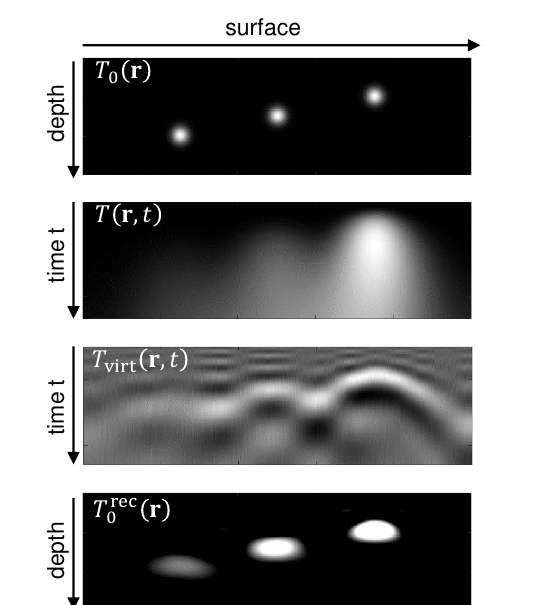}
    \caption{Two-dimensional reconstructions with simulated data using the virtual wave concept. $T_0(\mathbf{r})$ above shows the initial temperature distribution, $T(\mathbf{r},t)$ is the "measured" surface temperature at a depth of $z=0$. The surface temperature signal has a $SNR$ of approx. 1000. $T_\text{virt}(\mathbf{r},t)$ is calculated from the inversion of Eq. (\ref{Eq:Temperature_virtual_wave}) by the truncated SVD method, from which $T_0^\text{rec}(\mathbf{r})$ is calculated by applying the SAFT back-projection method for reconstruction from an acoustic wave.}
    \label{fig:VWC}
\end{figure}

For three-dimensional reconstruction, the $SNR$ cannot be improved by direct averaging of the surface signals. However, similar to reconstruction with ultrasound, reconstruction methods can be used that include the diffusion of heat in all directions. This provides a kind of averaging over all detection points and the $SNR$ can be increased similar to direct averaging. Acoustic imaging methods, such as back-projection, Synthetic Aperture Focusing Technique (SAFT), or time reversal, can be used by calculating in a first step from the measured temperature the so-called virtual wave, which has the same initial heat or pressure distribution as the measured temperature, but is a solution of the wave equation \cite{Burgholzer.2017}. The measured surface temperature and the virtual wave are connected by a local transformation, i.e. for the same location $\mathbf{r}$:
\begin{eqnarray}
    T(\mathbf{r},t) &&= \int_{- \infty}^\infty T_\text{virt}(\mathbf{r},t') M_T(t,t')\text{d}t', \nonumber \\
    && \text{with} \nonumber \\
    M_T(t,t') &&\equiv \frac{c}{\sqrt{\pi \alpha t}} \exp{\left( - \frac{c^2 t'^2}{4 \alpha t} \right)} \phantom{X} \text{for} \phantom{X} t>0 \,,
    \label{Eq:Temperature_virtual_wave}
\end{eqnarray}
where the sound velocity $c$ of the virtual wave can be chosen arbitrarily. Fig. \ref{fig:VWC} shows a two dimensional example with simulated data where the virtual wave concept is used for reconstruction of the original temperature distribution. Here, for the initial temperature distribution three gaussian-like temperature peaks are located at different depths below the sample surface.  The deepest source is barely visible when looking at the 200 measured surface temperature signals $T(\mathbf{r},t)$ alone with a SNR of 1000, whereas when reconstructed using the virtual wave concept, all three original temperature sources are clearly visible in $T_0^\text{rec}(\mathbf{r})$. The ultrasound reconstruction from the virtual waves has the same effect on the SNR as averaging and therefore, compared to a single one-dimensional reconstruction, in this example for 200 detection points the SNR can be enhanced be by a factor of $\sqrt{200}$ or the spatial resolution is enhanced by a factor of $\ln{\sqrt{200}} \approx 2.65$ according to Eq. (\ref{Eq:delta_resolution_thermal}). Since the $SNR$ increases with the square root of the averaged measurements and the resolution increases with $\ln(SNR)$ (Eq. (\ref{Eq:delta_resolution_thermal})), increasing the $SNR$ by averaging or by using ultrasound reconstruction from virtual waves, increases the resolution only in a logarithmic way. \par

\begin{figure}[htb]
\centering 
\includegraphics[trim={7.5cm 3.0cm 7.2cm 2.5cm},clip,width=0.6\textwidth]{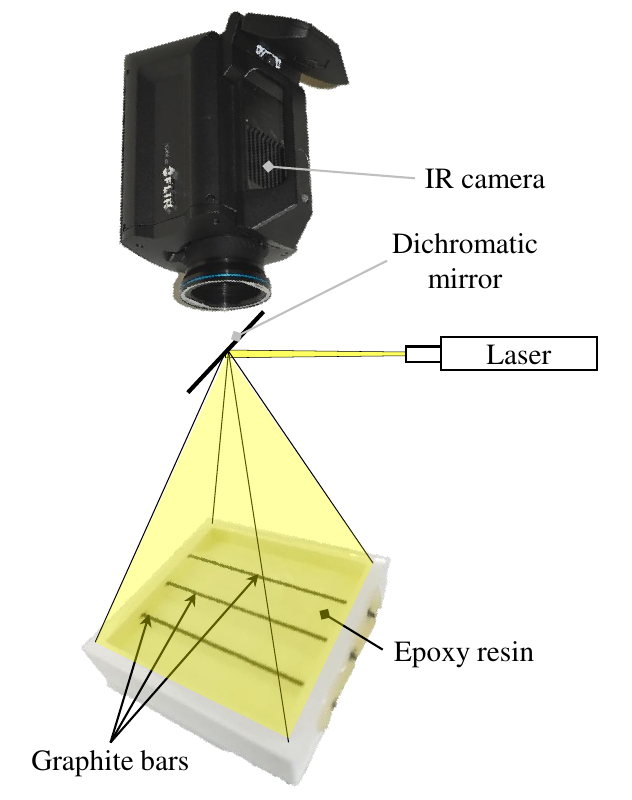}
    \caption{Sketched measurement set-up with three graphite bars embedded in epoxy resin at a depth of 1.6 mm, 2.6 mm, and 3.6mm, a laser for a short heating of the bars, and an infrared camera to measure the temporal evolution of the surface temperature. (Adapted from G. Thummerer et al., Photoacoustics 19, 100175, 2020; licensed under Creative Commons Attribution (CC BY) license.)}
    \label{fig:Bars}
\end{figure}

\begin{figure}[h]
\centering
\includegraphics[trim={11.2cm 9.0cm 1.3cm 9.0cm},clip,width=0.7\textwidth]{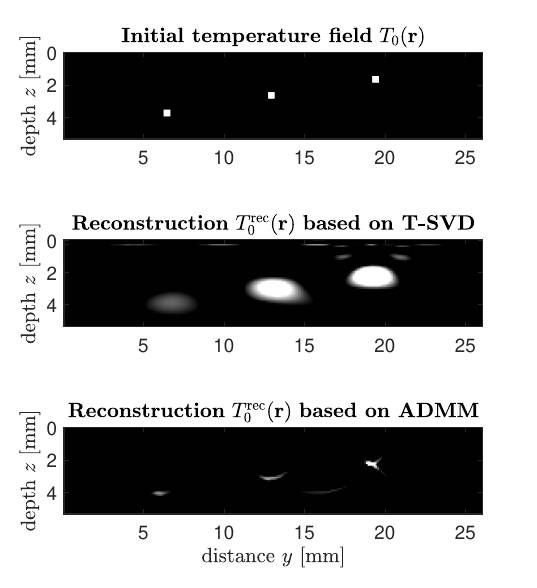}
    \caption{Comparison of the initial temperature distribution and the reconstructed initial temperature distributions using T-SVD or ADMM for calculating the virtual wave from the surface temperature measurement, followed by a subsequent reconstruction using SAFT. (Adapted from G. Thummerer et al., Photoacoustics 19, 100175, 2020; licensed under Creative Commons Attribution (CC BY) license.)}
    \label{fig:Barsrec}
\end{figure}

To increase the thermal resolution limit more efficiently, either super-resolution methods, such as those known to overcome the Abbe limit in optical imaging, e.g., by structured illumination, or additional information, such as sparsity or positivity, can be used. When using structured illumination, several reconstructed images with different illumination are used to calculate one super-resolution image by a non-linear reconstruction method, such as the iterative joint sparsity algorithm (IJOSP) \cite{Burgholzer.2017b}. As an example of the consideration of sparsity and positivity, the reconstruction of three graphite bars embedded in epoxy resin is presented here (thermal diffusivity $0.13 \times 10^{-6}$ m$^2$/s at a depth of 1.6 mm, 2.6 mm, and 3.6mm) \cite{Thummerer.2020b}. The measurement set-up is sketched in Fig. \ref{fig:Bars}. The graphite bars were heated by laser excitation (diode laser wavelength 938 nm, 250 W laser power with a pulse duration of 200 ms). The surface temperature after laser pulse excitation was measured with an infrared camera (106 Hz full frame mode and a noise equivalent temperature difference of 25 mK, cooled InSb sensor), which was sensitive in the spectral range of $3.0 - 5.1$ $\mu$m. In this spectral range the epoxy resin is opaque and we measured the surface temperature evolution. The spatial resolution from the camera pixels was approx. 0.1 mm on the sample surface. Fig. \ref{fig:Barsrec} shows a cross section of the initial temperature field from the light absorbing bars, and the reconstruction using the truncated-singular value decomposition (T-SVD) method for calculating the virtual wave and SAFT for the ultrasound reconstruction algorithm. To calculate the virtual wave considering sparsity (the cross section of graphite rods are only small points) and positivity (heating always increases the temperature), the Alternating Direction Method of Multipliers (ADMM) was used, a nonlinear optimization algorithm similar to the Douglas-Rachford method. Sparsity is taken into account by these methods by minimizing in addition to a data fitting term also the L1-norm.  More precisely, the virtual wave is computed as a minimizer  of $\| \mathbf{M} \, T_{\text{virt}} - T\|^2 + \mathcal{R}(T_{\text{virt}}) $, where $\mathbf{M}$ is a  discretization of the integral in  (\ref{Eq:Temperature_virtual_wave}) and $\mathcal{R}(.)$ is a suitable regularizer defined by the L1-norm. The use of this regularized virtual wave allows to implement efficient L1-minimization algorithms (e.g. ADMM) and leads to a much better resolution of the reconstructed temperature field \cite{Thummerer.2020b}.\par

The application of deep neural networks to tackle inverse problems is currently attracting a huge amount of interest within the inverse problems community, e.g. for thermal property reconstruction \cite{Glorieux.1996,Glorieux.1999,Glorieux.1999b} or for source reconstruction \cite{Ravi.2005}. Combining these methods with traditional model-based reconstruction techniques for solving inverse problems is an emerging trend which comes with new challenges, such as generating training data, designing architectures, choosing learning strategies, and interpreting the learned parameters. For recent results in this field, we refer to the survey of S.~Arridge et al.~\cite{arridge2019solving}. We have compared purely model based reconstructions using the ADMM method for calculating the virtual wave and a subsequent reconstruction with SAFT, a hybrid reconstruction by using deep learning to determine the reconstruction from the virtual wave, and a end-to-end deep learning reconstruction using directly the measured temperature data without calculating the virtual wave before \cite{SPM}. The best spatial resolution could be demonstrated by the hybrid reconstruction using the ADMM method to calculate the virtual wave and then the deep learning reconstruction. This showed better results than the reconstruction directly by deep learning and also the effort for training of the network was significantly less. For training of the deep neural network (u-net architecture) we have used simulated synthetic data, but nevertheless the reconstruction from measured data was better than using our model-based approach \cite{JAPdeeplearning}. One reason for this could be that the virtual waves show clearer signals with more spatial variations, similar to the images that are to be reconstructed, and not as slightly varying signals as the original thermal measurements. Using simulated synthetic data for training is important, because real-world data are hard to obtain, as it requires production of physical \textit{phantoms} with a range of material properties (e.g., defects at various positions and of various sizes and shapes).
However, we can synthesize arbitrarily large amounts of data with relatively simple analytic models of the heat diffusion process,
albeit at the risk of results based on synthetic data not perfectly matching the real world.
As a consequence, we provided code and data (synthetic and real-world) for full reproducibility of the results obtained by our experiments.\par

In summary, the virtual wave concept is an ideal and universal tool enabling the prospect of high-resolution photothermal imaging. The virtual wave concept can be viewed as beneficial factorization of the highly ill-posed inverse photothermal problem into a conversion from thermal to acoustic waves and a subsequent beamforming task. The conversion task is still highly ill-posed. However, the challenging situation of the full 2D or 3D reconstruction problem has been reduced to the solution of an 1D integral equation Eq. (\ref{Eq:Temperature_virtual_wave}). In doing so, we were able to efficiently integrate sparsity as well as positivity of the initial temperature profile to overcome the resolution limit. Note that the virtual waves $T_\text{virt}(\mathbf{r},t)$ are generally not positive. However, if we apply the Abel transform to the virtual waves, we obtain the spherical mean values of the initial temperature profile. Therefore, in reality, positivity is always satisfied for the Abel transformed waves that we actually use for the  positivity constraint. In the pixel basis, sparsity may also not be given directly. However, sparsity of the initial temperature is usually the case in various different bases such as the wavelet base. As shown in \cite{zangerl2021multiscale} in the context of photoacoustic imaging, wavelet sparsity of the initial temperature is transferred to wavelet sparsity of the virtual wave, thereby allowing the use of sparse recovery techniques. Finally, some other information such as the shape of absorbers may also be available. Exploiting such information is interesting and challenging. We anticipate, however, that ideas from \cite{zangerl2021multiscale} along with matched filtering techniques are promising for integrating such prior information into the reconstruction of the virtual wave.

\section{\label{sec:Discussion}Discussion}

Often, high frequency components in the signals cannot be detected because the noise prevents them from detection. There are several effects that limit the detectability and cause noise. Beside insufficient instrumentation and data processing one principle limitation comes from thermodynamic fluctuations (noise) associated with dissipative processes in accordance with the dissipation-fluctuation theorem (Fig. \ref{fig:coins}). As mentioned above, in addition to the imaging depth, the $SNR$ also affects spatial resolution and the width of the thermal PSF. Even if the noise amplitude is the same in frequency and time domain for Gaussian white noise, because the noise energy must be the same in the frequency and time domain (Parseval's theorem), the signal amplitude and hence the $SNR$ can be very different. For example, a short delta-shaped pulse in time contains all frequencies with the same amplitude and therefore the $SNR$ is much lower in the frequency domain than in the time domain. The appropriate space for the $SNR$ to calculate the resolution limit is closely related to the ill-posed inverse problem to compensate for the effect of thermal diffusion or acoustic attenuation. Regularization is required to solve such ill-posed inverse problems \cite{Hansen.1998,Aster.2018,Scherzer.2009}. The use of the cutoff frequency $\omega_\text{cut}$ is one way to regularize the ill-posed inverse problem, which turns out to be equivalent to the truncated singular value decomposition (SVD) method. But of course such a rigorous cut, which uses the full information of the frequency components below $\omega_\text{cut}$ and completely neglects them above this cutoff frequency, is only a rough approximation. In frequency domain, the information content of a frequency component gradually decreases with the square of the reduced amplitude as the frequency increases, until it disappears in the noise at the cutoff frequency $\omega_\text{cut}$. Therefore, other regularization methods that allow a smoother transition, such as Tikhonov regularization \cite{Burgholzer.2011}, might better describe this behavior.\par

To describe how humans perceive acoustic waves, Harvey Fletcher conducted tone-in-noise masking experiments, where a single tone (now a frequency domain peak equivalent to the delta time peak from above) was masked by noise. He had a problem with at the time (very) noisy telephones and needed to find out why people could or could not hear the speech and what the telephone company could do about it \cite{Fletcher}. We performed some of these experiments with different signals and with white noise and bandstop noise to better understand whether signals are detectable or not in the presence of simultaneous noise. It turned out that before the time-domain signal is converted to the frequency domain, temporal windowing is an important aspect in describing detectability. This in turn affects the possible measurement bandwidth for the signal, and a trade-off between higher bandwidth or more noise must be made in choosing the optimal discretization time step. This could lead to an optimal temporal discretization to achieve an optimal spatial resolution.\par

The physical reason for the entropy production are fluctuations, which is the noise added to an ideally noise-free signal and can be statistically described by stochastic processes \cite{Burgholzer.2013}. The ratio between the ideal signal amplitude and the noise amplitude is called the signal-to-noise ratio (SNR) and plays an essential role in this tutorial. Noise is identified as the thermodynamic fluctuations around a certain time varying mean-value and is the mechanism for entropy production and information loss. Other technical limitations from insufficient instrumentation can additionally limit spatial resolution more than described by the thermal PSF.\par

The thermodynamic fluctuations in macroscopic samples are usually so small that they can be neglected – but not for ill-posed inverse problems. The fluctuations are highly amplified due to the ill-posed problem of image reconstruction. As long as the macroscopic mean-value-equations describe the mean entropy production, the resolution limit depends only on the amplitude of the fluctuations and not on the actual stochastic process including all correlations \cite{Burgholzer.2015}. Therefore, we have already used in the past a simple Gauss-Markov process to describe the measured signal as a time-dependent random variable \cite{Burgholzer.2013}. Here, we describe heat diffusion as a random walk (Wiener process).\par

\section{\label{sec:Conclusions}Conclusions}
This tutorial explains the origin of the spatial resolution limit associated to photothermal techniques for imaging buried structures. The diffusion process is represented as a stochastic process (random walk) that causes spatial spreading of the mean values (entropy increase) and is accompanied by fluctuations that lead to information loss about the original temperature distribution, imposing a resolution limit that is calculated. Thereby it shows how the diffusion of heat and the fluctuations leads to the well-known blurring of reconstructions in infrared imaging. This can be described mathematically by a spatial convolution of the original temperature distribution with a thermal point-spread-function (PSF). The width of the thermal PSF increases linearly with imaging depth and is indirectly proportional to the natural logarithm of the signal-to-noise ratio (SNR).\par

The random walk model is introduced as a mathematically simple approach based on statistical probability calculations that provides insight into the heat diffusion process. It is also explained how probabilities come into play in a deterministic model for diffusion, where interaction with the environment is the cause for loosing information. Furthermore, some inversion strategies based on the virtual wave concept, the inclusion of sparsity and positivity of the solution, as well as deep neural networks, are presented and compared as suitable methods to partially overcome the depth-dependent resolution limit. 

\section*{Acknowledgments}
This work is co-financed by research subsidies granted by the government of Upper Austria. Financial support was also provided by the Austrian research funding association (FFG) under the scope of the COMET programme within the research project “Photonic Sensing for Smarter Processes (PSSP)” (contract number 871974). This programme is promoted by BMK, BMDW, the federal state of Upper Austria and the federal state of Styria, represented by SFG. Parts of this work have been supported by the Austrian Science
Fund (FWF), projects P 30747-N32 and P 33019-N. In addition, the financial support by the Austrian Federal Ministry for Digital and Economic
Affairs and the National Foundation for Research, Technology and Development and the Christian Doppler
Research Association is gratefully acknowledged.

\section*{Data availability}
The data that support the findings of this study are available from the corresponding author upon reasonable request.


\end{document}